\documentclass{article}

\usepackage{arxiv}

\usepackage[utf8]{inputenc} 
\usepackage[T1]{fontenc}    
\usepackage{hyperref}       
\usepackage{url}            
\usepackage{booktabs}       
\usepackage{amsfonts}       
\usepackage{nicefrac}       
\usepackage{microtype}      
\usepackage{lipsum}

\usepackage{graphicx}
\usepackage{subfigure}
\usepackage{amsmath}
\usepackage{tikz}
\usepackage{enumitem}
\usepackage{arydshln}
\usepackage{algorithm}
\usepackage{algorithmic}
\usepackage{hyperref}

\title{Solving All Regression Models For Learning Gaussian Networks\\
			Using Givens Rotations }

\author{
  Borzou Alipourfard\\
  Department of Computer Science\\
  Univesity of Texas at Arlington, Texas\\
  Arlington, USA \\
  \texttt{borzou.alipourfard@mavs.uta.edu} \\
   \And
  Jean X. Gao \\
  Department of Computer Science\\
  Univesity of Texas at Arlington, Texas\\
  Arlington, USA \\
  \texttt{gao@uta.edu} \\
}

\begin{document}
\maketitle

\begin{abstract}
Score based learning (SBL) is a promising approach for learning Bayesian networks. The initial step in the majority of the SBL algorithms consists of computing the scores of all possible child and parent-set combinations for the variables. For Bayesian networks with continuous variables, a particular score is usually calculated as a function of the regression of the child over the variables in the parent-set. The sheer number of regressions models to be solved necessitates the design of efficient numerical algorithms. In this paper, we propose an algorithm for an efficient and exact calculation of regressions for all child and parent-set combinations. In the proposed algorithm, we use QR decompositions (QRDs) to capture the dependencies between the regressions for different families and Givens rotations to efficiently traverse through the space of QRDs such that all the regression models are accounted for in the shortest path possible.  We compare the complexity of the suggested method with different algorithms, mainly those arising in all subset regression problems, and show that our algorithm has the smallest algorithmic complexity. We also explain how to parallelize the proposed method so as to decrease the runtime by a factor proportional to the number of processors utilized.\footnote{This work has been submitted to a journal for possible publication. Copyright may be transferred without notice, after which this version may no longer be accessible.}
\end{abstract}


\section{Introduction}
Bayesian networks (BN) are used to portray probabilistic dependencies. They are graphical models that encode a set of conditional independence statements through absence or presence of directed edges among nodes in a graph \cite{lauritzen1996graphical, jensen1996introduction}.  The child and parent relations formed in a such a graph help capture the dependency structure of the domain. BNs have found great application in machine learning, knowledge modeling, desicion systems, and causal learning \cite{spirtes1989causality, friedman2004inferring, pearl2009causality}.

A significant statistical problem is to learn a BN from observational data.  A promising approach for tackling this learning problem consists of a group of algorithms under the heading of score based learning algorithms (SBL) \cite{cooper1992bayesian, heckerman1995learning, chickering2002finding}. The initial step in the majority of the SBL algorithms consists of computing the \textit{local} scores for all possible child and parent-set combinations \cite{chickering2002optimal, silander2012simple, teyssier2012ordering}.  In the case of Bayesian networks over continous domain, the local score for a particular child and parent-set pair is usually calculated as a function of the regression of the child variable over the parent-set.  The sheer number of regressions models that need solving presents a computational challenge. For a network with $m$ nodes, there are $2^{m-1}$ candidate parent sets for a particular node and a total of $m2^{m-1}$ candidate families (a family consisting of a particular child and parent-set pair). The exponential number of node and parent-set combinations neccesiates the need for efficient numerical strategies.

 In this paper, we propose an algorithm for an efficient and exact calculation of regressions for all child and parent-set combinations of a given set of variables. The main data structure employed in the proposed algorithm is QR decomposition (QRD) as it both provides high numerical percision and can capture the dependency between the regressions for different families.  Using QRD together with Givens rotations (a low-cost operation), we show how to form a sequence of \textit{R} matrices that provide the necessary information required to solve the regressions for all the families.  In the proposed algorithm, we find the shortest of such sequences.

We further compare the theoretical runtime of the suggested method with different algorithms, mainly those arising in all subset regression problems, and show that our algorithm has the fastest runtime. We also explain how to parallelize the proposed algorithm providing a linear speedup proportional to the number of processors utilized.

\subsection{Previous Work}
A brute strategy for calculating all regression models would be to solve the regression equations for each family independent of the other families. Given the QRD of a family with $k$ parents, the number of flops (a flop consisting only of basic arithmetic operations such as a multipication, a division, a subtraction, or an addition) required for calculating the regression coefficients is of the order of $O(k^2)$. Forming the QRD itself requires $O(2nk^2)$ flops given $n$ samples. Therefore, the total number of flops required for calculating the regressions for all parent-set and child pairs is:

\begin{equation}
\label{eq:1}
O(naive) = m  \sum_{k=1}^{m} (2nk^2 + k^2)C(k,m), 
\end{equation}

where $C(k,m)$ is the number of possible combinations of $k$ objects out of $m$, and $m$ is the number of variables.

We can achieve a faster runtime by using the covariance matrix of the data and Cholesky decompositions and form each QRD in $O(k^3/3)$. This however comes at cost of numerical accuracy \cite{ortega1990numerical}. Using (\ref{eq:1}) and replacing $2nk^2 + k^2$ with $k^3/3$, it can be shown that forming the QRDs for all parent-set and child pairs requires $O(m^2 2^{m-2})$ flops. 

A better strategy would be to consider the problem of calculating the regression models for all the families as that of solving for multiple all-subset regression problems \cite{miller2002subset, smith1989all, yanev2004algorithms}. In other words, calculating the regression models for a BN over $m$ variables can be thought of as $m$ all-subset regression problems, one for each node. Using algorithms specialized to all-subset regression problems such as Clarke's algorithm or DCA, it is possible to achieve a linear speed up in solving for regression models \cite{gatu2007graph, gatu2003parallel, clarke1981algorithm}. However, such algorithms still solve the all-subset regression problem for each node independent of the other nodes and do not utilize the full structure of the problem. Thus we conjecture that it is possible to improve the runtime even further.

\subsection{Givens Rotations}
Reviewing the basics of Givens rotations would help us explain the workings of our proposed algorithm better. Pre-mutliplying a vector with a Givens rotation matrix, $G_\theta^{(i, j)}$, rotates the vector in the $(i,j)_{th}$ plane:

\begin{small}
\begin{gather}
G_\theta^{(i, j)} = \left(
\begin{array}{c:ccc:c}
I &  & 0 &  & 0 \\
\hdashline
 & \cos(\theta)  &  & \sin(\theta)  &  \\
 0 & & & & 0 \\
 & -\sin(\theta)  &  & \cos(\theta)  & \\
\hdashline
0 &  & 0 &  & I \\
\end{array} \right),
\end{gather}
\end{small}

where $\cos(\theta)$ and $\sin(\theta)$ in $G_\theta^{(i, j)}$ appear at the intersections of $i_{th}$ and $j_{th}$ rows and columns. By choosing $\theta$ appropriately, it is possible to rotate a vector such that its $j_{th}$ component becomes zero while preserving its norm [7]. 

\subsection{Retriangularization with Givens Rotation}
Consider a data matrix for the variables $\{ v_1, .., v_i, v_{i+1}$ $, ..., v_m\}$ and its corresponding QRD (with the same variable ordering). Having this QRD, it is then possible to compute the QRD of the data matrix with variables ordered as $\{ v_1, ..,  v_{i+1}, v_i, ..., v_m\}$ (variables $\{i \}$ and $\{i+1\}$ are transposed in the new order) by pre-multiplying the $R$ factor (after having its $\{i \}$ and $\{i+1\}$ columns swapped) by an appropriate Givens rotation matrix:

\begin{gather}
\hat{R} = G_\theta^{(i, i+1)} \left(
\begin{array}{cccccc}
r_1,1    & \ldots & r_{1,i+1}     & r_{1,i}  & \ldots & r_{1, m} \\
0          & \ldots & r_{2,i+1}     & r_{2,i}  & \ldots & r_{2, m} \\
\vdots  &           & \vdots        & \vdots  &           & \vdots \\
0          & \ldots & r_{i,i+1}     & r_{i,i}    & \ldots & r_{i, m} \\
0          & \ldots & r_{i+1,i+1} & 0          & \ldots & r_{i+1, m} \\
\vdots   &           & \vdots        & \vdots  &           & \vdots \\
0           & \ldots & 0               & 0          &  \ldots  & r_{m,m} \\
\end{array} \right), 
\end{gather} 

where $\theta$ can be calculated using :

\begin{equation}
\begin{array}{cc}
\cos(\theta) & = \dfrac{r_{i,i+1}}{\sqrt{(r_{i,i+1}^2 + r_{i+1,i+1}^2)}} \\
			&			\\
\sin(\theta) & = \dfrac{r_{i,i+1}}{\sqrt{(r_{i,i+1}^2 + r_{i+1,i+1}^2)}} \\
\end{array} 
\end{equation}

Such matrix transformation requires $O(6*(m-i+1))$ flops. We can regard this procedure as an operator that given a starting QRD outputs another QRD by column swapping and retriangularization. We will call this operator $GRC$ and will refer to it by $\xrightarrow{\text{GRC}}$. We use $\xrightarrow{\text{GRC}_i}$ to emphasize that column swapping occurs at index $i$. 

\section{Problem Definition}
\begin{figure}[t]
\centering
\begin{tikzpicture}
	[scale=1,auto=left,every node/.style={circle,fill=blue!0}]
	\node(a) at (-2,0) {$^{[A,B,C]}$};
	\node(b) at (0,0) {$^{[B,A,C]}$};
	\node(c) at (2,0) {$^{[B,C,A]}$};
	\node(d) at (2,-2) {$^{[C,B,A]}$};
	\node(e) at (0,-2) {$^{[C,A,B]}$};
	\node(f) at (-2,-2) {$^{[A,C,B]}$};
	
	\draw[thick, <->]  (a) -- (b) node[midway,above]  {$_{_{GRC}}$} ;
	\draw[thick, <->]  (b) -- (c) node[midway,above]  {$_{_{GRC}}$} ;
	\draw[thick, <->]  (c) -- (d) node[midway,right]  {$_{_{GRC}}$} ;
	\draw[thick, <->]  (d) -- (e) node[midway,below]  {$_{_{GRC}}$} ;
	\draw[thick, <->]  (e) -- (f) node[midway,below]  {$_{_{GRC}}$} ;
	\draw[thick, <->]  (f) -- (a) node[midway,left] {$_{_{GRC}}$} ;
	
\end{tikzpicture}
\caption{Each node represents a QRD. An edge indicates that one can calculate the QRD of a neighbouring node through using GRC operation.}
\label{Fig:1}
\end{figure}
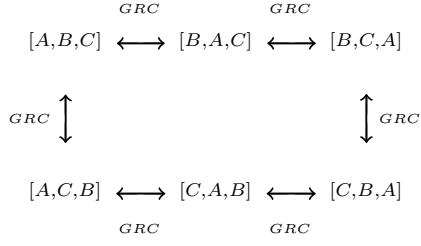

We will describe our problem's framework and our proposed algorithm by considering an example Bayesian network with three nodes, $\{A,B,C\}$. We use this example to construct a more general framework by the end of this section. In a BN over three variables, there are a total of nine regression models: $\{[A]\perp B, [A]\perp C, [B]\perp A, [B]\perp C, [C]\perp B,[C]\perp A, [A,B]\perp C, [A, C]\perp B, [B, C]\perp A \}$, where $[x,y]\perp z$ represents the regression model with $\{x, y\}$ as predictor variables and $z$ as the response variable.

For these variables, there are also six QRDs, each corresponding to a specific premutation of these variables. Each QRD can, in turn, be used to solve three regression models:

\begin{equation}
[A,B,C] \Rightarrow \{[A]\perp B, [A]\perp C, [A,B] \perp C \},
\end{equation}

where we have used $[A,B,C]$ for the QRD of the data matrix with variables ordered as $[A,B,C]$. We can go from one permutation to another (from one QRD to another) by using the operation introduced above and transposing two adjacent columns and retriangularization of the resulting matrix by using Givens transformation:

\begin{equation}
\left[ A,B,C \right]  \xrightarrow{\text{GRC}_2}  [A, C, B].
\end{equation}

The GRC operation allows us to traverse between the six QRDs for the variables $[A,B,C]$ as shown in the graph of Fig~\ref{Fig:1}.

Now consider the sequence of column transpositions giving rise to sequence of permutations:

\begin{gather}
\begin{array}{ccc}
[A, B, C]  & \xrightarrow{\text{GRC}_1} [B, A, C]  \xrightarrow{\text{GRC}_2} [B, C, A]  \\
 & \xrightarrow{\text{GRC}_1} [C, B, A]  \xrightarrow{\text{GRC}_2} [C, A, B]. 
\end{array}
\end{gather} 

Note that all 9 regression models are included in this permutation set. Therefore, it is possible to obtain all the QRDs necessary for solving all the regression models in a BN with three variables starting with an initial QRD, and forming four more by four adjacent column transpositions and application of four Givens rotations. Thus the problem of solving for all nine regression models is reduced to performing four simple Givens rotations.

Using this example as a basis, we propose the following algorithm. The algorithm starts with calculating the QRD of the data matrix for $m$ variables. This initial QRD can be used to solve regression models for $\dfrac{m(m-1)}{2}$ families. Then, through a sequence of adjacent column transpositions and retriangularizations, new QRDs are calculated. Each new QRD provides the information required to solve further regression models until all regression models are accounted for. The only remaining decision is to choose the sequence of column transpositions optimally such that minimum number of column transpositions are employed. 

\section{A Greedy Algorithm Using Givens Rotations} 
 In our algorithm, we choose the sequence of column transpositions in a greedy manner, leading to a very simple and intuitive algorithm. The greedy choice at each step consists of chooing a column for swapping such that the number of new models specified by the new permutation is the highest. Following is the sequence of greedy column transpositions that provides a sequence of QRDs sufficient for solving all the regressions of a BN with four variables, $ \{ A, B, C, D\}$ (the total number of models specified at each step is shown next to the permutations):
 
\begin{equation}
 \begin{array}{ccc}
\left[ A, B, C, D \right] & \Rightarrow & (6 \: Models) \\
\left[ B, A, C, D \right] & \Rightarrow & (9 \: Models)\\
\left[ B, C, A, D \right] & \Rightarrow & (11 \: Models)\\
\left[ C, B, A, D \right] & \Rightarrow & (14 \: Models)\\
\left[ C, A, B, D \right] & \Rightarrow & (16 \: Models)\\
\left[ C, A, D, B \right] & \Rightarrow & (17 \: Models)\\
\left[ C, D, A, B \right] & \Rightarrow & (19 \: Models)\\
\left[ D, C, A, B \right] & \Rightarrow & (22 \: Models)\\
\left[ D, A, C, B \right] & \Rightarrow & (24 \: Models)\\
\left[ D, A, B, C \right] & \Rightarrow & (25 \: Models)\\
\left[ D, B, A, C \right] & \Rightarrow & (26 \: Models)\\
\left[ D, B, C, A \right] & \Rightarrow & (28 \: Models)
\end{array}
\end{equation}

On closer examination, one will see that the above greedy algorithm follows a recursive structure. More specifically, consider the variables $ \{ X_1, X_2, ..., X_{m-1}, X_{m}\}$. Assume that we are given the QR decomposition of the variables in the same order as written above. Note that the greedy algorithm starts at the left most position. Further, note that transposing at $X_{m}$, the last variable, at any step, leads to a permutation that only adds one single model. Therefore, a greedy algorithm can always limit its operation to transposing the first $m-1$ variables until all the $m(2^{m-2}-1)$ regression models having predictors in $\{ X_1, X_2, ..., X_{m-1} \}$ are solved. When the permutations are such exhausted that transpositions in these positions add no new models, then the first transposition of $X_m$ occurs: $ \{ \hat{X_1}, \hat{X_2}, ..., X_{m}, \hat{X_{m-1}}\}$. In the following steps, unless $X_m$ is at first position, then no transposition on other variables is allowed since such transpositions add no new models. This sequence continues until through transpositions on $X_m$, this variable comes to the first position of the permutation, $ \{ X_m, \hat{X_1}, \hat{X_2}, ..., \hat{X_{m-1}} \}$. At this step, a greedy algorithm for $ \{ X_m, \hat{X_1}, \hat{X_2}, ..., \hat{X_{m-1}}\}$ applies the same sequence of transpositions to $\{\hat{X_1}, \hat{X_2}, ..., \hat{X_{m-1}}\}$ that it previously applied to $\{ X_1, X_2, ..., X_{n-1}\}$. Therefore, we can implement the proposed greedy algorithm through recursion. 

In particular, assume that we have found a sequence of greedy column transpositions that generates a set of QRDs sufficient for finding all the $m(2^{m-1} - 1)$ regression models for a BN with $m$ nodes. Let us call this sequence of column transpositions $\Upsilon(m)$. The greedy sequence of column transpositions for a graph with $m+1$ nodes can then be formed in three steps. (1)  $\Upsilon(m+1)$ starts with the same sequence of column transpositions as that of $\Upsilon(m)$, leading to permuatations necessary for calculation of the regressions of every node $X_i$ over all possible parent-sets not containing $X_{m+1}$. (2) To calculate the regressions of nodes $X_i$ over parent-sets containing $X_{m+1}$, we first move the variable $X_{m+1}$ to the start of the ordering by applying $m$ column transpositions. (3) The final sequence of column swapping in $\Upsilon(m+1)$ consists of swappings at indexes $i+1$ for $i \in \Upsilon(m)$. The resulting sequence of column transposition then form the sequence of greedy column transpositions that generate for finding permutaitons necessary the score table of a BN with $m+1$ nodes. The pseudocode GreedySwaps in Algorithm~\ref{alg:GreedySwaps} uses this recursive structure to find the sequence of greedy swaps.

\begin{algorithm}[tb]
\caption{GreedySwaps$(m)$}
\label{alg:GreedySwaps}
\begin{algorithmic} 
\STATE {\bfseries Input:} number of nodes in the BN $m$ 
\STATE {\bfseries Ouput:} An array of length $2^m -m -1$ of swapping indexes
\IF{$m == 2$}
\STATE {\bfseries return}  $[1]$
\ELSE
\STATE $\hat{\Upsilon} \leftarrow$ GreedySwaps$(m-1)$ 
\STATE {{\bfseries return}  $(\hat{\Upsilon}: [m-1..1] : [i | i \in \hat{\Upsilon}])$}
\STATE \#{ concatanation operator is shown by $:$}
\ENDIF
\end{algorithmic}
\end{algorithm}

Clarke has proposed an algorithm for solving all subset regression problem with a similar recursive structure \cite{clarke1981algorithm, smith1989all}. There are, however, major differences between the two algorithms. First, in every QRD, we consider \textit{all} possible combinations of regressor and predictor variables. For example, in the case of having the QRD of variables $\{A, B, C, D\}$, given in alphabetic ordering, Clarke is only concerned with the three regression models $\{[A, B, C] \bot D, [A, B] \bot D, [A] \bot D \}$, while we consider the regression models $\{[A, B] \bot C, [A] \bot B, [A] \bot C \}$ in addition to that of Clarke's. Furthermore, as discussed later in section 6, our algorithmic complexity analysis shows that the proposed algorithm results in a linear speed up compared to that of Clarke's.

Using the recurrence relation described above, we can find the length of the greedy sequence of column transpositions as a function of number of variables in the BN:

\begin{equation}
\begin{array}{ccc}
|\Upsilon(m)| = & 2|\Upsilon(m-1)| + m-1, \\
|\Upsilon(m)| = & 2^m - 1 - m.
\end{array}
\end{equation}

\section{Optimality of Greedy Algorithm}
In this section, we show that the proposed greedy algorithm is optimal; that the number of $GRC$ operations or column transpositions required for generating a sufficient set of QRDs for solving all the regression models of a BN using the greedy algorithm is minimal. Our proof makes use of an auxiliary problem. This problem is that of solving all subset regression problem for $m$ predictors and one (imaginary) response variable. We first show that our original problem is equivalent to this auxiliary problem when the only operation allowed is that of column transposition and retriangularization. Since the optimal solution for the auxiliary problem is known [6], we can then show that our solution is optimal for the original problem.

In the first step of our proof we show that the two following problems are equivalent:

\begin{enumerate}[label=(\roman*)]
\item Problem I (Original). Forming a sequence of QRDs such that all the regression models for a BN with $m$ variables are included in the QRD set. The constraints are that we only have access to the starting QRD with variables ordered as $[X_1, ..., X_m]$. Further, the only operation available is transposing two immediate column and retriangularization using the $GRC$ operator.

\item Problem II (Auxiliary). Forming a sequence of QRDs such that all the regression models in an all subset rgression problem with $m$ predictor variables are included in the QRD set. Again we are given a starting QRD with variables ordered as $[X_1, ..., X_m, Y]$ ($Y$ is the response variable). We are also again constrained to only using the $GRC$ operator.
\end{enumerate}

We first show that every solution to Problem I is also a solution to problem II. Specifically, assume that we are looking to solve the regression of $Y$ over predictor set $P = \{X_{1P}, .., X_{KP}\}$. In other words, we want a QRD where the variables, $V= \{X_{1}, .., X_{m}\}$, are ordered as $[perm_1(P), perm_2(V\backslash P)]$, and $perm_1(P)$ and $perm_2(V \backslash P)$ are some permutations of variables in the sets $P$ and $V \backslash P$ . Assume that the size of the predictor set $P$ is smaller than $m$. Since we know the solution to Problem I, then we have access to QRDs that provide the solution to the regressions of node $X_j$, $X_j \in V\backslash P$, over all its possible parent sets. Thus, the solution to Problem I, has to provide a QRD of the form $[\hat{perm_1}(P), \hat{perm_2}(V \backslash P)]$. Choosing $perm_1(P)$ and $perm_2(V \backslash P)$ equal to $\hat{perm_1}(P)$ and $\hat{perm_2}(V \backslash P)$, we get the desired QRD. If $P == V$, we simply choose a node $X_j$, and the solution to Problem I, provides access to a QRD of the form $[perm_1(V \backslash \{X_j\}), X_j]$, which is the desired QRD for solving problem II.

In the same manner, we can show that every solution to Problem II is also a solution to Problem I. Assume we wish to solve a regression for a specific parent set and node pair where the parents are in the $F = \{X_{1F}, .., X_{KF}\}$, and the node is $X_j \notin F$. In other words, we want a QRD where the variables are ordered as $[perm_1(F), perm_2(V \backslash F)]$, where $V= \{X_{1}, .., X_{m}\}$, and $perm_1(F)$ and $perm_2(V \backslash F)$ are some permutations of variables in the sets $F$ and $V \backslash F$. Given the solution to Problem II, we have QRDs for all subset regression models. Therefore, we have access to a QRD where the first $m$ variables are ordered as $[\hat{perm_1}(F), \hat{perm_2}(V\backslash F)]$. Choosing $perm_1(P)$ and $perm_2(V \backslash F)$ equal to $\hat{perm_1}(P)$ and $\hat{perm_2}(V \backslash F)$, we get the desired QRD. This concludes our proof that the above two problems are equivalent.

Our solution to Problem I, uses $2^{m} - m -1$  column transpositions. It has been proven that the minimal number of column transpositions required for solving Problem II, is in fact $2^{m} - m -1$. Therfore, the proposed greedy algorithm is optimal and no other algorithm can generate a sequence of QRDs of smaller length such that all regression models of a score table for a BN are included in the QRD set. 

\section{Algorithmic Complexity Analysis}
In this section we calculate the runtime of our algorithm and compare it to three other methods for solving the regression models for all possible families of a BN.  To analyze the runtime of the proposed algorithm, we make use of the following recurrence relation for the sequence of column transpositions employed by our algorithm: 

\begin{equation}
\begin{array}{cc}
\Upsilon(m) =  & \Upsilon(m-1): [m-1, ..., 1]  \\
 & :[i+1\mid i \in \Upsilon(m-1)], 
\end{array}
\end{equation}

where $\Upsilon(m)$ is the sequence of column transpositions for $m$ variable case, $[m-1..1]$ is a sequence of column transpositions starting at index $m-1$ down to the first position, and $:$ is concatanation operator. Noting that the number of operations required for applying Given's transformation to a matrix of size $m$ is six more than its counterpart for a matrix of size $m-1$, we can write the following recurrence to describe the runtime of our algorithm:

\begin{equation}
\begin{array}{cc}
\label{eq:5}
T(m) = &[T(m-1)+6|\Upsilon(m-1)|]   \\
&+ \sum_{i=1}^{m-1}6(m-i+1)+ T(m-1),
\end{array}
\end{equation}

where $|\Upsilon(i)|$ is the number of column transpositions employed by our algorithm when applied to a network with $i$ nodes:

\begin{equation} 
\label{eq:6}
\begin{aligned}
|\Upsilon(i)| = 2^i-1-i.
\end{aligned}
\end{equation}

Combining these two equations, (\ref{eq:5}) and (\ref{eq:6}), we can find the runtime of our algorithm to be of the order:

\begin{equation}
\begin{aligned}
T(m) = 3m2^m + 62^m-3m^2-9m-6.
\end{aligned}
\end{equation}

The runtime of the proposed method for forming the QRDs of all possible parent-set and child combination is compared to that of Clakre's all subset regression algorithm \cite{clarke1981algorithm}, Dropping Column Algorithm (DCA) \cite{gatu2003parallel}, and direct brute force using Cholesky decomposition and covariance structure in Table~\ref{table_1}.

\begin{table}[t]
\caption{Comparison of the runtime of the proposed method (in bold) to three other algorithms.}
\label{table_1}
\vskip 0.15in
\begin{center}
\begin{small}
\begin{sc}
\begin{tabular}{|c||c|}
\hline
Algorithm & Runtime \\
\hline
\hline
\textbf{Greedy Column Swapping} & $\boldsymbol{O(3m2^m)}$\\
\hline
DCA & $O(9m2^m)$\\
\hline
Clarke & $O(1.5m^22^m)$ \\
\hline
Brute Force & $O(0.5m^32^m)$\\
\hline
\end{tabular}
\end{sc}
\end{small}
\end{center}
\vskip -0.1in
\end{table}

\section{Parallelization}
In order to design an efficient parallel version of the proposed algorithm, let us re-state the general problem framework introduced in section 4 so as to account for employing of multiple processors. Given $p$ processors, we wish to find $p$ initial QRDs and $p$ sequences of adjacent column transpositions for each of the processors, such that the resulting variable permutations and their corresponding QRDs among all the processing cores solves for all the $m*(2^{m-1} - 1)$ regression models of a Bayesian network with $m$ nodes. 

To find the optimal performance gain using $P$ processors, we first derive a bound on possible performance improvment in the case of all subset regression problem. Due to discussions in section 5, we know that this bound would be still in effect for the problem of forming all QRDs for all families in a BN. We then propose a near optimal parallelization scheme that achieves a performance gain close to this bound.

Assume that we have found $p$ initial QRDs, $R_i$, $i =1, ..., p$ and $p$ sequences of column transpositions, $\Upsilon_i$, $i =1, ..., p$, for solving the all subset regression problem using $p$ processors.  We will denote the length of the sequence of column transpositions performed by processor $i$ by $|\Upsilon_i|$. Since the $p$ initial QRDs can at most account for $mp$ regression models and since a column transposition can at most add one new regression model, we have:

\begin{equation}
\sum_{i=1}^{p} |\Upsilon_i| > 2^m - 1 - mp.
\end{equation}

Therefore we have:
\begin{equation}
\exists i: |\Upsilon_i| > \dfrac{2^m - 1 - mp}{p}.
\end{equation}

Thus the number of required column transpositions when using $p$ processors is at best of the order of $O(\dfrac{2^m}{p})$.

We propose a parallelization method that is a direct consequence of the following recurrence relation:

\begin{equation}
\begin{array}{cc}
\Upsilon_0(m) =  & \Upsilon_0(m-1): [m-1, ..., 1]  \\
 & :[i+1\mid i \in \Upsilon_0(m-1)], 
\end{array}
\end{equation}

where $\Upsilon_0(m)$ denotes the sequence of column transpositions performed by the single core algorithm of the previous section for the $m$ variable case.

Consider the case where the number of processors, $p$, is two. In this case, we initialize the first core with the QRD of variables ordered as ${v_1, ..., v_m}$ and the second core with the QRD of variables ordered as ${v_m, v_1 ..., v_{m-1}}$. For the first processor, we choose the sequence of column transpositions equal to $\Upsilon_0(m-1)$, and for the second processor we employ the sequence of column transpositions $[i+1\mid i \in \Upsilon_0(m-1)]$. In general, given $p=2^k$ processors, we can use this recurrence relation $k$ times to find the initial QRDs and the sequence of column transpositions for each processor. 

More specifically let us represent each of the $2^k$ processors with a binary array of $k$-bits, mapping processor $i$ to a binary array equivalent to its binary repressentation. Then, Algorithm~\ref{alg:par} can provide the initial permutation and the sequence of column swapping required to be performed at processor $i$.

\begin{algorithm}[tb]
   \caption{SeedPathCalculator$(P_{id}, m)$}
   \label{alg:par}
\begin{algorithmic}
   \STATE {\bfseries Input:} processor binary array code $P_{id}$, size of the Bayesian network $m$, 
   \STATE {\bfseries Output:} initial permutation $\upsilon$, sequence of column swapping indexes $\Upsilon$
   \STATE $k = m-len(P_{id})$
   \STATE $\upsilon = [1..k]$
   \FOR{$i=k+1$ {\bfseries to} $n$} 
		\IF{$P_{id}[m-i] == 0$}
			\STATE $\upsilon.insert(0, i)$
		\ELSE
			\STATE $\upsilon.append(i)$
		\ENDIF
	\ENDFOR
	\STATE $d = sum(P_{id}) $
	\STATE $\Upsilon = d + GreedySwaps(k)$
	\STATE {\bfseries return} $(\upsilon, \Upsilon)$
\end{algorithmic}
\end{algorithm}

The number of column transpositions performed by each processors in general is of the order of $O(2^{m-k}-1-(m-k))$ or $O(\dfrac{2^m}{p})$. 

\section{Conclusion}
In this paper, we proposed an algorithm for an efficient and exact calculation of regressions for all the families of a BN. Noting that the regressions for the different families are dependent on each other, we utilized QR decomposition as a data structure for capturing these dependencies. We then used Givens rotations and column transpositions as low-cost operations to efficiently trace a greedy path through the space of QRDs such that all the regression models are included. We showed how the proposed greedy method could be more easily implemented using recursions in section 3. In section 4 we provided a lower bound on the number of column transpositions required for solving the regressions for all the families and showed that the proposed greedy algorithm achieves this lower bound. 

We further compared the runtime of our algorithm with specialized algorithms for all-subset regression problems in Table~\ref{table_1}. We argued that spcialized all-subset regression algorithms and brute force algorithms do not utilize the whole of the dependency structure among the families. The faster runtime of our proposed method then proves that we make better use of the dependency structure. Although in terms of algorithmic complexity our proposed algorithm has only a constant factor of improvement compared to that of DCA, the memmory requirements of our algorithm is much lower than theirs. Specifically, the proposed algorithm utilizes a storage of size $O(m^2)$ (only a single $R$ matrix needs to be stored at any moment) while DCA requires a storage of size $O(2^{m-3})$. In section 7 we further provided a near optimal parallelization scheme for our prosed algorithm. 


\bibliography{/Users/Borzou/Documents/Research/My_Publications/Arxiv_Given/Arxiv_Given_Bilb}

\begin{thebibliography}{10}

\bibitem{lauritzen1996graphical}
Steffen~L Lauritzen.
\newblock {\em Graphical models}, volume~17.
\newblock Clarendon Press, 1996.

\bibitem{jensen1996introduction}
Finn~V Jensen.
\newblock {\em An introduction to Bayesian networks}, volume 210.
\newblock UCL press London, 1996.

\bibitem{spirtes1989causality}
Peter Spirtes, Clark~N Glymour, and Richard Scheines.
\newblock {\em Causality from probability}, volume 112.
\newblock Carnegie-Mellon University, Laboratory for Computational Linguistics,
  1989.

\bibitem{friedman2004inferring}
Nir Friedman.
\newblock Inferring cellular networks using probabilistic graphical models.
\newblock {\em Science}, 303(5659):799--805, 2004.

\bibitem{pearl2009causality}
Judea Pearl.
\newblock {\em Causality}.
\newblock Cambridge university press, 2009.

\bibitem{cooper1992bayesian}
Gregory~F Cooper and Edward Herskovits.
\newblock A bayesian method for the induction of probabilistic networks from
  data.
\newblock {\em Machine learning}, 9(4):309--347, 1992.

\bibitem{heckerman1995learning}
David Heckerman, Dan Geiger, and David~M Chickering.
\newblock Learning bayesian networks: The combination of knowledge and
  statistical data.
\newblock {\em Machine learning}, 20(3):197--243, 1995.

\bibitem{chickering2002finding}
David~Maxwell Chickering and Christopher Meek.
\newblock Finding optimal bayesian networks.
\newblock In {\em Proceedings of the Eighteenth conference on Uncertainty in
  artificial intelligence}, pages 94--102. Morgan Kaufmann Publishers Inc.,
  2002.

\bibitem{chickering2002optimal}
David~Maxwell Chickering.
\newblock Optimal structure identification with greedy search.
\newblock {\em Journal of machine learning research}, 3(Nov):507--554, 2002.

\bibitem{silander2012simple}
Tomi Silander and Petri Myllymaki.
\newblock A simple approach for finding the globally optimal bayesian network
  structure.
\newblock {\em arXiv preprint arXiv:1206.6875}, 2012.

\bibitem{teyssier2012ordering}
Marc Teyssier and Daphne Koller.
\newblock Ordering-based search: A simple and effective algorithm for learning
  bayesian networks.
\newblock {\em arXiv preprint arXiv:1207.1429}, 2012.

\bibitem{ortega1990numerical}
James~M Ortega.
\newblock {\em Numerical analysis: a second course}.
\newblock SIAM, 1990.

\bibitem{miller2002subset}
Alan Miller.
\newblock {\em Subset selection in regression}.
\newblock Chapman and Hall/CRC, 2002.

\bibitem{smith1989all}
DM~Smith and JM~Bremner.
\newblock All possible subset regressions using the qr decomposition.
\newblock {\em Computational Statistics \& Data Analysis}, 7(3):217--235, 1989.

\bibitem{yanev2004algorithms}
Petko Yanev, Paolo Foschi, and Erricos~John Kontoghiorghes.
\newblock Algorithms for computing the qr decomposition of a set of matrices
  with common columns.
\newblock {\em Algorithmica}, 39(1):83--93, 2004.

\bibitem{gatu2007graph}
Cristian Gatu, Petko~I Yanev, and Erricos~J Kontoghiorghes.
\newblock A graph approach to generate all possible regression submodels.
\newblock {\em Computational Statistics \& Data Analysis}, 52(2):799--815,
  2007.

\bibitem{gatu2003parallel}
Cristian Gatu and Erricos~J Kontoghiorghes.
\newblock Parallel algorithms for computing all possible subset regression
  models using the qr decomposition.
\newblock {\em Parallel Computing}, 29(4):505--521, 2003.

\bibitem{clarke1981algorithm}
MRB Clarke.
\newblock Algorithm as 163: A givens algorithm for moving from one linear model
  to another without going back to the data.
\newblock {\em Journal of the Royal Statistical Society. Series C (Applied
  Statistics)}, 30(2):198--203, 1981.

\end{thebibliography}
\bibliographystyle{unsrt}

\end{document}